# High-performance computation of large sparse matrix exponential

Feng Wu *, Kailing Zhang, Li Zhu, and Jiayao Hu

*State Key Laboratory of Structural Analysis of Industrial Equipment,*

*Department of Engineering Mechanics, Faculty of Vehicle Engineering and*

*Mechanics, Dalian University of Technology, Dalian 116023, P.R.China*

*Email: wufeng_chn@163.com, zkl1008@mail.dlut.edu.cn, zhuli@mail.dlut.edu.cn,*

*HJYdxyx@mail.dlut.edu.cn*

*Corresponding author: Tel: 8613940846142; e-mail address: wufeng_chn@163.com

**Abstract:** Computation of the large sparse matrix exponential has been an important topic in many fields, such as network and finite-element analysis. The existing scaling and squaring algorithm (SSA) is not suitable for the computation of the large sparse matrix exponential as it requires greater memories and computational cost than is actually needed. By introducing two novel concepts, i.e., real bandwidth and $\varepsilon$- bandwidth, to measure the sparsity of the matrix, the sparsity of the matrix exponential is analyzed. It is found that for every matrix computed in the squaring phase of the SSA, a corresponding sparse approximate matrix exists. To obtain the sparse approximate matrix, a new filtering technique in terms of forward error analysis is proposed. Combining the filtering technique with the idea of keeping track of the incremental part, a competitive algorithm is developed for the large sparse matrix exponential. The proposed method can primarily alleviate the over-scaling problem due to the filtering technique. Three sets of numerical experiments, including one large matrix with a dimension larger than $2\times10^6$, are conducted. The numerical experiments show that, compared with the expm function in MATLAB, the proposed algorithm can provide higher accuracy at lower computational cost and with less memory.

# 1. Introduction

For $\boldsymbol{H}\in\mathbb{C}^{n\times n}$, the matrix exponential $e^{\boldsymbol{H}}$ can be defined as

$$e^{\boldsymbol{H}}=\sum_{i=0}^{+\infty}\frac{1}{i!}\boldsymbol{H}^{i}\ , \tag{1.1}$$

which is a matrix function widely used in many engineering fields, such as, control theory [1], structural vibration [2], heat conduction [2], networks[4], image recognition[5] and quantum physics[6], and is an area of strong focus in academia. To calculate the matrix exponential accurately and efficiently, many methods have been put forward, such as Taylor series expansion, Padé approximation, and the scaling and squaring algorithm (SSA)[7]. In Ref. [8], 20 methods for calculating the matrix exponential are systematically summarized and compared.

When the dimension of the matrix exponential to be calculated is large (e.g., greater than 100), the SSA method is often preferred because it has the advantages of high precision and good efficiency and thus attracts the attention of many scholars. Many improved SSAs have been proposed in recent years. In Ref. [9], the over-scaling problem of the SSA was alleviated by using $\left\|\boldsymbol{H}^k\right\|^{1/k} \leq \|\boldsymbol{H}\|$. In Refs. [10-13], the SSA was combined with new algorithms based on the Taylor polynomial to improve the computational efficiency. In Ref. [14], an accurate and efficient SSA was developed with matrix splines. In Ref. [15], a fast algorithm combining the SSA with low displacement rank was proposed for the Toeplitz matrix exponential. In Ref. [16], a sharp bound was placed on the SSA and the number of matrix multiplication was reduced, especially matrix multiplications in the squaring phase. In Ref. [16] a new calculation formula for calculating the squaring process was proposed that can reduce the rounding error caused by the finite arithmetic in the squaring process, which is called the precise integration method (PIM), and it is widely used in structural vibration analysis[2, 17]. The improvement of the SSA method mentioned above mainly focus on error analysis of the SSA method and reducing the number of matrix multiplications but do not involve large-scale matrix exponential calculations. However, in many fields, such as structural vibration analysis, heat transfer, and networks, a matrix exponential calculation of a large sparse matrix is often involved, the matrix dimension of which is often greater than 10,000. For the large-scale matrix, more attentions are focused on the product of $e^{\boldsymbol{H}}\boldsymbol{v}$ [18, 19], and the studies on the calculation of $e^{\boldsymbol{H}}$ for a large-scale matrix are much less. One main reason is that the matrix exponential $e^{\boldsymbol{H}}$ is often dense, even $\boldsymbol{H}$ is sparse. However, $e^{\boldsymbol{H}}$ usually contains more physical information than $e^{\boldsymbol{H}}\boldsymbol{v}$ and hence it is very valuable to obtain $e^{\boldsymbol{H}}$ in some applications [1, 4, 20, 21]. Ref. [4] pointed out that $e^{\boldsymbol{H}}$, instead of $e^{\boldsymbol{H}}\boldsymbol{v}$, is needed to evaluate the communicability of complex network, and it is of practical significance to compute $e^{\boldsymbol{H}}$ by utilizing the sparsity of matrix. The elements in $e^{\boldsymbol{H}}$ often exhibit certain decay properties. Some decay rates of the elements in $e^{\boldsymbol{H}}$ have been established in Refs. [22-25]. In Refs. [3, 26], the authors explain the decay properties based on the physical properties that wave speed is finite, and a fast precise integration method (FPIM) was proposed. In the squaring phase, the FPIM will discard the elements with extremely small absolute values in the matrix, so as to obtain a sparser matrix and significantly improve computational efficiency (in this paper, this technique is called filtering). The FPIM is a promising method that can efficiently calculate the large-scale matrix exponential involved in vibration response and non-Fourier heat conduction, and realize the exponential calculation of matrices with tens of thousands of dimensions. However, the following questions still remain to be answered with respect to this method

1） How will the sparsity of the matrix change during the squaring phase? Answering this question helps to explain why the filtering technique is useful.

2） How to identify which elements in the matrix can be discarded and which elements should be retained?

3） When considering filtering techniques, how to determine the number of squarings and the number of Taylor series expansion terms?

In this paper, the above questions are discussed in detail, and an efficient and high-precision algorithm for calculating the large-scale sparse matrix exponential is presented. The rest of this paper is organized as follows. In Section 2, we first review the FPIM. In Section 3, we focus on discussing the sparsity of the matrix, introduce the concepts of real bandwidth and $\varepsilon$-bandwidth, and discuss several properties of aforementioned concepts. On this basis, we study the sparsity of the matrix exponential, laying a mathematical foundation for the filtering technology. In Section 4, the forward error analysis of the PIM based on a Taylor series is presented, and on this basis the filtering technology is discussed. In combination with the sparsity of the matrix exponential, the selection strategy of the number of Taylor series expansion terms and the number of squares considering the effect of the filtering technology is given. In Section 5, we give some numerical examples of large sparse matrix exponentials, among which the dimension of the largest matrix exceeds $2\times10^6$. The analyses in Sections 3 and 4 are verified through numerical experiments. Finally, conclusions are given in Section 6.

For all matrices $\boldsymbol{A}$ in this paper, we always assume $\boldsymbol{A}\in\mathbb{C}^{n\times n}$, where $n$ represents the dimension of the matrix, and its matrix elements are denoted $a_{ij}$. The norm $\|\boldsymbol{A}\|$ adopted in this paper refers to any compatible norm unless otherwise stated.

# 2. Existing PIMs

The PIM was first proposed to solve structural vibration analysis problems, e.g.,

$$\dot{\boldsymbol{U}}=\boldsymbol{\Theta}\boldsymbol{U},\quad \boldsymbol{U}=\begin{pmatrix}\boldsymbol{a}\\ \boldsymbol{v}\end{pmatrix},\quad \boldsymbol{\Theta}=\begin{bmatrix}\boldsymbol{0} & \boldsymbol{I}\\ -\boldsymbol{M}^{-1}\boldsymbol{K} & -\boldsymbol{M}^{-1}\boldsymbol{C}\end{bmatrix}, \tag{2.1}$$

where, $\boldsymbol{M}$, $\boldsymbol{C}$ and $\boldsymbol{K}$ are the lumped-mass matrix, damping matrix and stiffness matrix,

respectively, and $\boldsymbol{a}$ and $\boldsymbol{v}$ are the displacement and velocity vectors, respectively, and $\boldsymbol{I}$ is the unite matrix. It should be noted that the lumped-mass matrix is mandatory, otherwise $\boldsymbol{H}$ will be dense, which is not considered in this investigation. If the solution at time $t_1$ is denoted $\boldsymbol{U}_1$, the solution at time $t_2 = t_1 + \Delta t$ can be expressed as:

$$\boldsymbol{U}_2 = e^{\boldsymbol{H}}\boldsymbol{U}_1, \quad \boldsymbol{H} = \boldsymbol{\Theta}\Delta t . \tag{2.2}$$

It can be seen from Eq. (2.1) that the calculation of $e^{\boldsymbol{H}}$ plays an important role in solving Eq. (2.1). The easiest way to calculate the matrix exponential is the Taylor series expansion. However, according to Theorem 2.2 proposed in Ref. [16], if $\|\boldsymbol{H}\|$ is large, the rounding error of the Taylor series computation may be very serious. A specific example is found in Ref. [8]. The rounding error difficulty can be controlled by the SSA, which utilizes the following properties:

$$e^{\boldsymbol{H}} = \left[ e^{\boldsymbol{H}2^{-N}} \right]^{2^N} , \tag{2.3}$$

Letting $\boldsymbol{F}_i := e^{\boldsymbol{H}2^{(i-N)}}$, then obviously

$$\boldsymbol{F}_i = \boldsymbol{F}_{i-1} \cdot \boldsymbol{F}_{i-1}, \; e^{\boldsymbol{H}} = \boldsymbol{F}_N , \tag{2.4}$$

When $N$ is large enough, $\|\boldsymbol{H}2^{-N}\| \ll 1$, so it is easy to calculate $\boldsymbol{F}_0 = e^{\boldsymbol{H}2^{-N}}$ by a Taylor series or Padé approximants without serious rounding errors. When a Taylor series is used to calculate $\boldsymbol{F}_0$, the result of the SSA method can be expressed as

$$\boldsymbol{F}_i \approx F_{MNi}\left(\boldsymbol{H}\right) = \left( \sum_{s=0}^{M} \frac{1}{s!} \left(\boldsymbol{H}2^{-N}\right)^s \right)^{2^i} . \tag{2.5}$$

In terms of Eq. (2.5), $e^{\boldsymbol{H}} \approx F_{MNN}\left(\boldsymbol{H}\right)$, and

$$\boldsymbol{F}_0 \approx F_{M00}\left(\boldsymbol{H}\right) = \boldsymbol{I} + \boldsymbol{T}_0, \quad \boldsymbol{T}_0 = \sum_{s=1}^{M} \frac{1}{s!} \left(\boldsymbol{H}2^{-N}\right)^s . \tag{2.6}$$

However, if $N$ is large enough that the rounding error of the Taylor series computation is negligible, another type of rounding error may be unacceptable. Because $\|\boldsymbol{H}2^{-N}\| \ll 1$, $\|\boldsymbol{T}_0\| \ll 1$, and $\|\boldsymbol{I}\| = O(1)$. Adding large numbers to very small numbers will lead to serious rounding errors in the subsequent squaring process[27]. To reduce the rounding error, Zhong proposed the PIM [2] based on a key idea of "keeping track of the incremental part". It decomposes $F_{MNi}\left(\boldsymbol{H}\right)$ into $F_{MNi}(\boldsymbol{H}) = \boldsymbol{I} + \boldsymbol{T}_i$, the substitution of which into Eq. (2.4) yields:

$$\boldsymbol{T}_i = 2\boldsymbol{T}_{i-1} + \boldsymbol{T}_{i-1} \cdot \boldsymbol{T}_{i-1} \tag{2.7}$$

Carrying out the iterative steps in (2.7) $N$ times yields $\boldsymbol{T}_N$. Once $\boldsymbol{T}_N$ is obtained, $e^{\boldsymbol{H}} \approx \boldsymbol{I} + \boldsymbol{T}_N$. Compared with $\boldsymbol{I}$, $\boldsymbol{T}_N$ is no longer small, and therefore, there is no serious rounding error. The PIM is essentially an improved version of SSA. Compared with the SSA, the PIM adds one time of matrix addition in every iterative step, as it needs to add $2\boldsymbol{T}_i$ after calculating $\boldsymbol{T}_{i-1} \cdot \boldsymbol{T}_{i-1}$. Noting that the computation cost of the matrix addition is far smaller than that of matrix multiplication, the computation efficiency of PIM is very close to that of SSA, but the rounding error is reduced[2]. It should be pointed out that although the PIM was first proposed for structural vibration problems, it is also applicable to the high-precision calculation of the matrix exponential of other types of problems [2, 3, 17, 28].

Numerous practical calculation experiences using the PIM show that, although the PIM can obtain the matrix exponential $e^{\boldsymbol{H}}$ with high precision, it is still not suitable for large matrices (e.g., the dimension of matrix $\boldsymbol{H}$ is greater than 10,000). Aiming at the matrix exponential presented in Eq.(2.2), a FPIM is proposed in Refs. [3, 26]. For each matrix $\boldsymbol{T}_i$, the FPIM removes all the elements with very small absolute values based on an given filter parameter $\varepsilon_{\mathrm{f}}$. According to the numerical experiments, the computational efficiency of the FPIM is often improved by 2 orders of magnitude compared with the original PIM. However, the three questions on the FPIM, introduced in Section 1, still remain to be answered.

# 3. Sparsity of matrix exponentials

To answer the three questions, it is necessary to understand the sparsity of the matrix exponential. To analyze the sparsity of the matrix exponential, we first present two indicators, i.e., real bandwidth and $\varepsilon$- bandwidth, to measure the sparsity of a sparse matrix, and then discuss the sparsity of the matrix exponential using these two indicators.

## 3.1 Bandwidth of sparse matrix

For a sparse matrix $\boldsymbol{A}\in\mathbb{C}^{n\times n}$, the element of which is denoted $a_{ij}$, the upper and lower bandwidths $l_1(\boldsymbol{A})$ and $l_2(\boldsymbol{A})$ are defined by

$$l_1(\boldsymbol{A})=\max_{\Omega(\boldsymbol{A})}(j-i),\quad \text{and}\quad l_2(\boldsymbol{A})=\max_{\Omega(\boldsymbol{A})}(i-j), \tag{3.1}$$

respectively, where $\Omega(\boldsymbol{A}):=\{(i,j):\ a_{ij}\neq 0\}$. The bandwidth $l(\boldsymbol{A})$ of matrix $\boldsymbol{A}$ is defined by $l_1(\boldsymbol{A})+l_2(\boldsymbol{A})$. In this paper, it is always assumed that $l(\boldsymbol{A})\ll n$ for a sparse matrix.

**Lemma 3.1** *Let* $\boldsymbol{A}\in\mathbb{C}^{n\times n}$ *and* $\boldsymbol{B}\in\mathbb{C}^{n\times n}$ *be two sparse matrices, and* $\boldsymbol{C}=\boldsymbol{A}+\boldsymbol{B}$. *Then,*

(a) $l(\boldsymbol{A}+\boldsymbol{B})\leq \max\left(l_1(\boldsymbol{A}),l_1(\boldsymbol{B})\right)+\max\left(l_2(\boldsymbol{A}),l_2(\boldsymbol{B})\right)$;

(b) $l(\boldsymbol{AB})\leq l(\boldsymbol{A})+l(\boldsymbol{B})$.

***Proof*.** (a) Because

$$\begin{aligned}\Omega(\boldsymbol{A}+\boldsymbol{B})&=\{(i,j):\ a_{ij}+b_{ij}\neq 0\}\subset\{(i,j):\ |a_{ij}|+|b_{ij}|\neq 0\}\\&=\{(i,j):\ |a_{ij}|\neq 0\}\cup\{(i,j):\ |b_{ij}|\neq 0\}\\&=\{(i,j):\ a_{ij}\neq 0\}\cup\{(i,j):\ b_{ij}\neq 0\}=\Omega(\boldsymbol{A})\cup\Omega(\boldsymbol{B})\end{aligned},$$

we have

$$l_1(\boldsymbol{A}+\boldsymbol{B})=\max_{\Omega(\boldsymbol{A}+\boldsymbol{B})}(j-i)\leq\max_{\Omega(\boldsymbol{A})\cup\Omega(\boldsymbol{B})}(j-i)=\max\left(l_1(\boldsymbol{A}),l_1(\boldsymbol{B})\right)$$

and

$$l_2(\boldsymbol{A}+\boldsymbol{B})=\max_{\Omega(\boldsymbol{A}+\boldsymbol{B})}(i-j)\leq\max_{\Omega(\boldsymbol{A})\cup\Omega(\boldsymbol{B})}(i-j)=\max\left(l_2(\boldsymbol{A}),l_2(\boldsymbol{B})\right)$$

Hence, $l(\boldsymbol{A}+\boldsymbol{B})=l_1(\boldsymbol{A}+\boldsymbol{B})+l_2(\boldsymbol{A}+\boldsymbol{B})\leq\max\left(l_1(\boldsymbol{A}),l_1(\boldsymbol{B})\right)+\max\left(l_2(\boldsymbol{A}),l_2(\boldsymbol{B})\right)$.

(b) Let $\boldsymbol{C}=\boldsymbol{AB}$. According to Eq. (3.1), $a_{ik}=0$ when $k<i-l_2(\boldsymbol{A})$ or $k>i+l_1(\boldsymbol{A})$, and $b_{kj}=0$ when $k<j-l_1(\boldsymbol{B})$ or $k>j+l_2(\boldsymbol{B})$. As $c_{ij}=\sum_{k=1}^{n}a_{ik}b_{kj}$, $c_{ij}=0$ must hold when

$$j-l_1(\boldsymbol{B})>i+l_1(\boldsymbol{A}),\ \text{or,}\ j+l_2(\boldsymbol{B})<i-l_2(\boldsymbol{A}),$$

which means that a set consisting of $c_{ij}$, where $(i,j)$ satisfies

$$j-i\leq l_1(\boldsymbol{A})+l_1(\boldsymbol{B}),\ \text{or,}\ i-j\leq l_2(\boldsymbol{A})+l_2(\boldsymbol{B})$$

must contain all the nonzero elements in $C$. Hence,

$$l(C)=l_1(C)+l_2(C)\le l_1(A)+l_2(A)+l_2(B)+l_1(B)=l(A)+l(B). \qquad \blacksquare$$

We next introduce two parameters to describe the sparsity of the sparse matrix, i.e., the real bandwidth and the $\varepsilon$- bandwidth, defined as follows.

**Definition 3.2** *For a sparse matrix* $A\in\mathbb{C}^{n\times n}$, *and any* $\varepsilon>0$

(a) *the real bandwidth* $\lambda(A)$ *is*

$$\lambda(A)=\min_{P\in\mathcal{P}} l\left(PAP^{\mathrm{T}}\right), \tag{3.2}$$

*where* $\mathcal{P}$ *is the set of all the permutation matrices;*

*(b) the* $\varepsilon$*-bandwidth* $\beta(\varepsilon,A)$ *is*

$$\beta(\varepsilon,A)=\min_{B\in\Gamma(\varepsilon)}\lambda(A-B),\quad \Gamma(\varepsilon):=\left\{B:\ B\in\mathbb{C}^{n\times n},\ \text{and}\ \|B\|\le\varepsilon\|A\|\right\}. \tag{3.3}$$

**Lemma 3.3** For a sparse matrix $A\in\mathbb{C}^{n\times n}$, we have

(a) $\beta(\varepsilon,A)\le\lambda(A)\le l(A)$*;*

(b) *If* $\varepsilon_1\ge\varepsilon_2$, *then* $\beta(\varepsilon_1,A)\le\beta(\varepsilon_2,A)$.

***Proof****.* These two statements can be easily obtained from Definition 3.2. ■

Obviously, the real bandwidth and the $\varepsilon$-bandwidth are two indicators to measure the sparsity of the sparse matrix, and we have the following lemmas.

**Lemma 3.4** $A\in\mathbb{C}^{n\times n}$ *is a sparse matrix, and* $\bar{A}\in\mathbb{C}^{n\times n}$ *satisfies* $\|A-\bar{A}\|\le\delta\|A\|$*; then for any* $\varepsilon\ge\delta$,

$$\beta(\varepsilon,A)\le\beta(\bar{\varepsilon},\bar{A}),\quad \bar{\varepsilon}=\frac{\|A\|}{\|\bar{A}\|}(\varepsilon-\delta). \tag{3.4}$$

***Proof****.* According to Definition 3.2(b), a matrix $\bar{B}$ satisfying $\beta(\bar{\varepsilon},\bar{A})=\lambda(\bar{A}-\bar{B})$ and $\|\bar{B}\|\le\bar{\varepsilon}\|\bar{A}\|$ must exist; hence,

$$\beta(\bar{\varepsilon},\bar{A})=\lambda(\bar{A}-\bar{B})=\lambda\left(A-\left[(A-\bar{A})+\bar{B}\right]\right)\ge\beta\left(\frac{\left\|(A-\bar{A})+\bar{B}\right\|}{\|A\|},A\right).$$

Noting that $\frac{\left\|(A-\bar{A})+\bar{B}\right\|}{\|A\|}\le\frac{\|A-\bar{A}\|}{\|A\|}+\frac{\|\bar{B}\|}{\|A\|}\le\delta+\frac{\|\bar{B}\|}{\|\bar{A}\|}\frac{\|\bar{A}\|}{\|A\|}\le\delta+\bar{\varepsilon}\frac{\|\bar{A}\|}{\|A\|}$ and Lemma 3.3 (b), we have

$$\beta\left(\bar{\varepsilon},\bar{A}\right)\geq\beta\left(\frac{\left\|\left(A-\bar{A}\right)+\bar{B}\right\|}{\left\|A\right\|},A\right)\geq\beta\left(\delta+\frac{\left\|\bar{A}\right\|}{\left\|A\right\|}\bar{\varepsilon}\right)=\beta\left(\varepsilon,A\right)$$ ■

**Lemma 3.5** *Let* $A, B\in\mathbb{C}^{n\times n}$ *be sparse matrices, and* $p_m\left(A\right)=\sum_{i=0}^{m}c_iA^i$ *,* $c_i\in\mathbb{C}$*, where* $m$ *is an arbitrary positive integer. Then,*

(a) $l\left(p_m\left(A\right)\right)\leq ml\left(A\right)$*;*

(b) $\lambda\left(p_m\left(A\right)\right)\leq m\lambda\left(A\right)$*;*

(c) $\beta\left(\varepsilon_p,p_m\left(A\right)\right)\leq m\lambda\left(A-B\right)$, where $\varepsilon_p=\left\|p_m\left(A\right)\right\|^{-1}\left\|p_m\left(A\right)-p_m\left(A-B\right)\right\|$.

***Proof***. (a) Using Horner's method, the sum of $p_m\left(A\right)$ can be written as

$$B_1=c_{m-1}I+c_mA,\ B_i=c_{m-i}I+AB_{i-1},\quad p_m\left(A\right)=B_m=c_0I+AB_{m-1}.$$

Hence, we just need to prove that $l\left(B_m\right)\leq ml\left(A\right)$.

According to Lemma 3.1, and considering that $l\left(I\right)=0$, we have $l\left(B_1\right)\leq l\left(A\right)$. Supposing that $l\left(B_{i-1}\right)\leq\left(i-1\right)l\left(A\right)$ holds for $i-1$, and considering Lemma 3.1, we have $l\left(B_i\right)\leq l\left(AB_{i-1}\right)\leq l\left(A\right)+l\left(B_{i-1}\right)\leq i\times l\left(A\right)$. When $i=m$, $l\left(B_m\right)\leq ml\left(A\right)$.

(b) Let $P$ be the matrix satisfying $\lambda\left(A\right)=l\left(P^{\mathrm{T}}AP\right)$, then form Lemma 3.5(a),

$$\lambda\left(p_m\left(A\right)\right)\leq l\left(P^{\mathrm{T}}p_m\left(A\right)P\right)=l\left(\sum_{i=0}^{m}c_i\left(P^{\mathrm{T}}AP\right)^i\right)\leq ml\left(P^{\mathrm{T}}AP\right)=m\lambda\left(A\right).$$

(c) Letting $B_p=p_m\left(A\right)-p_m\left(A-B\right)$, then $\lambda\left(p_m\left(A\right)-B_f\right)=\lambda\left(p_m\left(A-B\right)\right)$. From Lemma 3.5(b), $\lambda\left(p_m\left(A-B\right)\right)\leq m\lambda\left(A-B\right)$. Noting that $\left\|B_p\right\|=\varepsilon_p\left\|p_m\left(A\right)\right\|$,

$$\beta\left(\varepsilon_p,p_m\left(A\right)\right)\leq\lambda\left(p_m\left(A\right)-B_p\right)\leq m\lambda\left(A-B\right).$$ ■

From Lemma 3.5(c), we can easily deduce the following.

**Corollary 3.6** *Let* $A\in\mathbb{C}^{n\times n}$ *be a sparse matrix, and* $p_m\left(A\right)=\sum_{i=0}^{m}c_iA^i$ *,* $c_i\in\mathbb{C}$*. Supposing* $B$ *is such that* $\beta\left(\varepsilon,A\right)=\lambda\left(A-B\right)$*, then*

$$\beta\left(\varepsilon_p,p_m\left(A\right)\right)\leq m\beta\left(\varepsilon,A\right),\quad\varepsilon_p=\frac{\left\|p_m\left(A\right)-p_m\left(A-B\right)\right\|}{\left\|p_m\left(A\right)\right\|}$$

Lemma 3.5 provides a bound for the bandwidth of an $m$ -order matrix polynomial. It is obvious that the bound increases linearly with $m$. As the matrix exponential $e^{\boldsymbol{H}}$ is an infinite series, the bound of the real bandwidth will be $n$, i.e., $e^{\boldsymbol{H}}$ may be a dense matrix, even if $\boldsymbol{H}$ is a sparse matrix. According to the SSA, the matrix exponential can be approximated precisely by a matrix polynomial, which means the $\varepsilon$- bandwidth of the matrix exponential may be small for some sparse matrices. Next, we discuss the $\varepsilon$- bandwidth of matrix exponential.

## 3.2 $\varepsilon$-bandwidth of matrix exponential

**Lemma 3.7** *Letting* $\boldsymbol{H} \in \mathbb{C}^{n\times n}$ *be a sparse matrix,* $\bar{\boldsymbol{H}}$ *be an approximate matrix of* $\boldsymbol{H}$*, and*

$$\varepsilon_{MNi}\left(\boldsymbol{H},\bar{\boldsymbol{H}}\right)=\frac{\left\|F_{MNi}\left(\boldsymbol{H}\right)-F_{MNi}\left(\bar{\boldsymbol{H}}\right)\right\|}{\left\|e^{\boldsymbol{H}2^{i-N}}\right\|}+\frac{\left\|e^{\boldsymbol{H}2^{i-N}}-F_{MNi}\left(\boldsymbol{H}\right)\right\|}{\left\|e^{\boldsymbol{H}2^{i-N}}\right\|}, \tag{3.5}$$

*where,* $F_{MNi}\left(\boldsymbol{H}\right)$ *is defined by Eq.* (2.5)*, then*

$$\beta\left(\varepsilon_{MNi}\left(\boldsymbol{H},\bar{\boldsymbol{H}}\right),e^{\boldsymbol{H}2^{i-N}}\right)\le M2^{i}\lambda\left(\bar{\boldsymbol{H}}\right). \tag{3.6}$$

***Proof.*** From Lemma 3.4, we have $\beta\left(\varepsilon_{MNi}\left(\boldsymbol{H},\bar{\boldsymbol{H}}\right),e^{\boldsymbol{H}2^{i-N}}\right)\le\beta\left(\bar{\bar{\varepsilon}}_{MNi},F_{MNi}\left(\boldsymbol{H}\right)\right)$, where

$$\bar{\bar{\varepsilon}}_{MNi}=\frac{\left\|e^{\boldsymbol{H}2^{i-N}}\right\|}{\left\|F_{MNi}\left(\boldsymbol{H}\right)\right\|}\left(\varepsilon_{MNi}\left(\boldsymbol{H},\bar{\boldsymbol{H}}\right)-\frac{\left\|e^{\boldsymbol{H}2^{i-N}}-F_{MNi}\left(\boldsymbol{H}\right)\right\|}{\left\|e^{\boldsymbol{H}2^{i-N}}\right\|}\right)=\frac{\left\|F_{MNi}\left(\boldsymbol{H}\right)-F_{MNi}\left(\bar{\boldsymbol{H}}\right)\right\|}{\left\|F_{MNi}\left(\boldsymbol{H}\right)\right\|}.$$

According to Lemma 3.5, $\beta\left(\bar{\bar{\varepsilon}}_{MNi},F_{MNi}\left(\boldsymbol{H}\right)\right)\le M2^{i}\lambda\left(\bar{\boldsymbol{H}}\right)$, therefore,

$$\beta\left(\varepsilon_{MNi}\left(\boldsymbol{H},\bar{\boldsymbol{H}}\right),e^{\boldsymbol{H}2^{i-N}}\right)\le\beta\left(\bar{\bar{\varepsilon}}_{MNi},F_{MNi}\left(\boldsymbol{H}\right)\right)\le M2^{i}\lambda\left(\bar{\boldsymbol{H}}\right). \quad \blacksquare$$

According to the above lemmas, the following theorem can be obtained.

**Theory 3.8** *Letting* $\boldsymbol{H} \in \mathbb{C}^{n\times n}$ *be a sparse matrix,* $\forall\varepsilon>0$*,* $\bar{\boldsymbol{H}}$ *be an approximate sparse matrix of* $\boldsymbol{H}$ *satisfying* $\left\|e^{\boldsymbol{H}2^{i-N}}-e^{\bar{\boldsymbol{H}}2^{i-N}}\right\|<\varepsilon\left\|e^{\boldsymbol{H}2^{i-N}}\right\|$*, and* $\varepsilon_{M_1N_1N_1}\left(\boldsymbol{H}2^{i-N},\bar{\boldsymbol{H}}2^{i-N}\right)$ *be defined by Eq. (3.5), then,*

$$\beta\left(\varepsilon,e^{\boldsymbol{H}2^{i-N}}\right)\le\bar{\alpha}\left(\boldsymbol{H}2^{i-N},\bar{\boldsymbol{H}}2^{i-N},\varepsilon\right)\lambda\left(\bar{\boldsymbol{H}}\right), \tag{3.7}$$

*where* $\bar{\alpha}\left(\boldsymbol{H}2^{i-N},\bar{\boldsymbol{H}}2^{i-N},\varepsilon\right):=\min\left\{M_1 2^{N_1}:\ \varepsilon_{M_1N_1N_1}\left(\boldsymbol{H}2^{i-N},\bar{\boldsymbol{H}}2^{i-N}\right)\le\varepsilon\right\}$.

***Proof.*** When $i=N$, inequality (3.7) becomes $\beta\left(\varepsilon,e^{\boldsymbol{H}}\right)\le\bar{\alpha}\left(\boldsymbol{H},\bar{\boldsymbol{H}},\varepsilon\right)\lambda\left(\bar{\boldsymbol{H}}\right)$, where $\bar{\alpha}\left(\boldsymbol{H},\bar{\boldsymbol{H}},\varepsilon\right):=\min\left\{M_1 2^{N_1}:\varepsilon_{M_1N_1N_1}\left(\boldsymbol{H},\bar{\boldsymbol{H}}\right)\le\varepsilon\right\}$. Noting that $\left\|e^{\boldsymbol{H}}-e^{\bar{\boldsymbol{H}}}\right\|<\varepsilon\left\|e^{\boldsymbol{H}}\right\|$, $\bar{\alpha}\left(\boldsymbol{H},\bar{\boldsymbol{H}},\varepsilon\right)$

always exists. According to Lemma 3.7, we have

$$\beta\left(\bar{\varepsilon}_{M_1N_1N_1}\left(\boldsymbol{H},\bar{\boldsymbol{H}}\right),e^{\boldsymbol{H}}\right)\le M_1 2^{N_1}\lambda\left(\bar{\boldsymbol{H}}\right),$$

which means $\beta\left(\varepsilon,e^{\boldsymbol{H}}\right)\le M_1 2^{N_1}\lambda\left(\boldsymbol{H}\right)$ holds for any $\left(M_1,N_1\right)$ satisfying $\bar{\varepsilon}_{M_1N_1N_1}\left(\boldsymbol{H},\bar{\boldsymbol{H}}\right)\le\varepsilon$. Hence, inequality (3.7) is true for $i=N$.

Noting that $\beta\left(\varepsilon,e^{\boldsymbol{H}}\right)\le\bar{\alpha}\left(\boldsymbol{H},\bar{\boldsymbol{H}},\varepsilon\right)\lambda\left(\bar{\boldsymbol{H}}\right)$ holds for any $\boldsymbol{H}$ and $\bar{\boldsymbol{H}}$ satisfying $\left\|e^{\boldsymbol{H}}-e^{\bar{\boldsymbol{H}}}\right\|<\varepsilon\left\|e^{\boldsymbol{H}}\right\|$, replacing $\boldsymbol{H}$ and $\bar{\boldsymbol{H}}$ with $\boldsymbol{H}2^{i-N}$ and $\bar{\boldsymbol{H}}2^{i-N}$, respectively, completes the proof. ■

According to Theorem 3.8, if $\bar{\alpha}\lambda\left(\bar{\boldsymbol{H}}\right)$ is far smaller than $n$, then there exists a sparse approximation matrix of $e^{\boldsymbol{H}2^{i-N}}$ with the relative error no greater than $\varepsilon$. We can properly filter $F_{MNi}\left(\boldsymbol{H}\right)$ to obtain a sparse matrix $\overline{F_{MNi}\left(\boldsymbol{H}\right)}$, so that the relative error of which is no greater than $\varepsilon$. Then, subsequent calculations based on $\overline{F_{MNi}\left(\boldsymbol{H}\right)}$ will require much less computational cost and memory. Theory 3.8 explains the rationality and necessity of a filtering process during the computation of the matrix exponential. How to filter is discussed in the next section.

# 4. Proposed method

Two issues should be considered about the filtering above. First, the effect on calculation accuracy should be controlled so that the numerical results should be correct; secondly, the filtered matrix should be as sparse as possible so that the best computational efficiency can be achieved. Both issues are detailedly discussed in this section.

## 4.1 Forward error upper bound of PIM

The PIM consists of two calculation phases, i.e., the calculation of $\boldsymbol{T}_0$ and the squaring processes. We first discuss the error in the squaring process. In the discussion, we ignore the effect of rounding error. For the convenience of discussion, Let $F_{MNi}\left(\boldsymbol{H}\right)$ be denoted as $\bar{\boldsymbol{F}}_i$, then $\bar{\boldsymbol{F}}_i=\bar{\boldsymbol{F}}_{i-1}^2=\bar{\boldsymbol{F}}_0^{2^i}$. Let the error matrix of $\bar{\boldsymbol{F}}_i$ be $\boldsymbol{R}_i=\boldsymbol{F}_i-\bar{\boldsymbol{F}}_i$, and considering that $\boldsymbol{R}_i$ and $\boldsymbol{F}_i$

commute, we have:

$$\begin{aligned} \boldsymbol{R}_i &= \boldsymbol{F}_i - \overline{\boldsymbol{F}}_i = \boldsymbol{F}_0^{2^i} - \overline{\boldsymbol{F}}_0^{2^i} = \boldsymbol{F}_0^{2^i} - \left(\boldsymbol{F}_0 - \boldsymbol{R}_0\right)^{2^i} \\ &= \boldsymbol{F}_0^{2^i} - \boldsymbol{F}_0^{2^i}\left(\boldsymbol{I} - \boldsymbol{F}_0^{-1}\boldsymbol{R}_0\right)^{2^i} = \boldsymbol{F}_0^{2^i}\left[\boldsymbol{I} - \left(\boldsymbol{I} - \boldsymbol{F}_0^{-1}\boldsymbol{R}_0\right)^{2^i}\right] \end{aligned} . \tag{4.1}$$

Then the relative error of $\overline{\boldsymbol{F}}_i$ will be

$$\begin{aligned} \frac{\left\|\boldsymbol{R}_i\right\|}{\left\|\boldsymbol{F}_0^{2^i}\right\|} &\le \left\|\boldsymbol{I} - \left(\boldsymbol{I} - \boldsymbol{F}_0^{-1}\boldsymbol{R}_0\right)^{2^i}\right\| = \left\|\boldsymbol{I} - \sum_{s=0}^{2^i}\binom{2^i}{s}(-1)^s\left(\boldsymbol{F}_0^{-1}\boldsymbol{R}_0\right)^s\right\| \\ &= \left\|-\sum_{s=1}^{2^i}\binom{2^i}{s}(-1)^s\left(\boldsymbol{F}_0^{-1}\boldsymbol{R}_0\right)^s\right\| \le \sum_{s=1}^{2^i}\binom{2^i}{s}\left\|\boldsymbol{F}_0^{-1}\boldsymbol{R}_0\right\|^s = \left(1 + \left\|\boldsymbol{F}_0^{-1}\boldsymbol{R}_0\right\|\right)^{2^i} - 1 = r_i \end{aligned} , \tag{4.2}$$

Formula (4.2) provides an upper bound $r_i$ for the relative error of $\overline{\boldsymbol{F}}_i$, and $r_i$ can actually be calculated through the iterative processes of the PIM, namely

$$r_0 = \left\|\boldsymbol{F}_0^{-1}\boldsymbol{R}_0\right\|, \quad r_i = 2r_{i-1} + r_{i-1}^2 . \tag{4.3}$$

Assuming that $r_{i-1}$ is smaller than the error tolerance, $r_{i-1}^2 \ll 2r_{i-1}$, the calculation of $r_i = 2r_{i-1} = 2^i r_0 = 2^i\left\|\boldsymbol{F}_0^{-1}\boldsymbol{R}_0\right\|$ is adopted in this paper. Then, the issue becomes the estimation of $r_0$. If $\boldsymbol{F}_0$ is computed by using the Taylor series, according to the formula (22) of [29], we have

$$\boldsymbol{F}_0^{-1}\boldsymbol{R}_0 = \sum_{i=0}^{+\infty}\left(\boldsymbol{H}2^{-N}\right)^{M+1+i}\frac{(-1)^i}{i!M!(i+M+1)} . \tag{4.4}$$

Hence, the norm of $\boldsymbol{F}_0^{-1}\boldsymbol{R}_0$ has the following bound:

$$\left\|\boldsymbol{F}_0^{-1}\boldsymbol{R}_0\right\| \le \frac{1}{M!}\sum_{i=0}^{+\infty}\frac{\left\|\boldsymbol{H}2^{-N}\right\|^{M+1+i}}{i!(i+M+1)} = \frac{(-1)^{M+1}}{M!}\gamma\left(M+1, -\left\|\boldsymbol{H}2^{-N}\right\|\right), \tag{4.5}$$

where $\gamma\left(M+1, -\left\|\boldsymbol{H}2^{-N}\right\|\right)$ is the lower incomplete gamma function.

In terms of formulas (4.3)-(4.5), the relative error upper bound of the PIM can be rewritten as

$$r_i = 2^i\frac{(-1)^{M+1}}{M!}\gamma\left(M+1, -\left\|\boldsymbol{H}2^{-N}\right\|\right) . \tag{4.6}$$

## 4.2 Filtering based on error analysis

As an approximation to $e^{\boldsymbol{H}2^{i-N}}$, $F_{MNi}(\boldsymbol{H})$ is a $M2^i$- order matrix polynomial, hence the upper bound of its real bandwidth is $M2^i\lambda(\boldsymbol{H})$, according to Lemma 3.5. Because $F_{M_1N_1N_1}\left(\boldsymbol{H}2^{i-N}\right)$ is also an approximation to $e^{\boldsymbol{H}2^{i-N}}$, according to Theorem 3.8, we have

$$\beta\left(r_i, e^{\boldsymbol{H}2^{i-N}}\right) \le \hat{\alpha}_i \lambda(\boldsymbol{H}), \tag{4.7}$$

where,

$$\hat{\alpha}_i := \min\left\{M_1 2^{N_1} :\ 2^{N_1}\frac{(-1)^{M_1+1}}{M_1!}\gamma\left(M_1+1, -\left\|\boldsymbol{H}2^{i-N}2^{-N_1}\right\|\right) \le r_i\right\}. \tag{4.8}$$

In the derivation from (4.4) to (4.5), we use the inequality $\left\|\boldsymbol{H}2^{-N(M+1+i)}\right\| \le \left\|\boldsymbol{H}2^{-N}\right\|^{M+1+i}$ which is over-scaling. In terms of Theorem 2 in [11], $\left\|\boldsymbol{H}2^{-N(M+1+i)}\right\| \le h_{M,N}^{M+1+i}$, where $h_{M,N} = \max_{M+1\le k\le 2M+1}\left\|\left(\boldsymbol{H}2^{-N}\right)^k\right\|^{\frac{1}{k}}$. Hence, $\hat{\alpha}_i$ can be re-evaluated as

$$\hat{\alpha}_i := \min\left\{M_1 2^{N_1} :\ 2^{N_1}\frac{(-1)^{M_1+1}}{M_1!}\gamma\left(M_1+1, -h_{M_1,N+N_1-i}\right) \le r_i\right\}. \tag{4.9}$$

It can be found from $\hat{\alpha}_i \le M2^i$. Actually, when $i \ne 0$, $\hat{\alpha}_i$ is generally much smaller than $M2^i$. A specific example is given here. We select 12 different normal matrices with the norm values $h \in [0,500]$. For each value, $(M,N)$ is the solution satisfying the condition $r_N \le 10^{-16}$ with the minimum of $M+N$, which is the times of matrix multiplication. $M2^i/\hat{\alpha}_i$ with different $i$ values are compared in Figure 1. It can be seen from Figure 1 that $\hat{\alpha}_i$ grows much more slowly than the $M2^i$. Hence for each $F_{MNi}(\boldsymbol{H})$ during the squaring phase, there exists a sparser approximate matrix with a relative error no greater than $r_i$, which demonstrates the necessity of filtering.

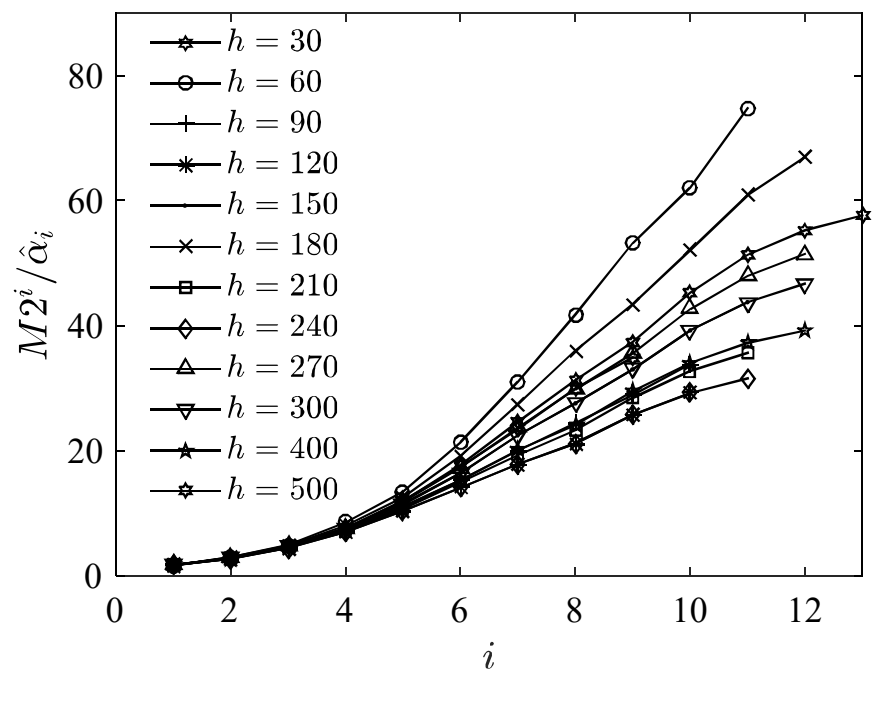

**Figure** 1 $M2^i/\hat{\alpha}_i$ for each norm

Since $\boldsymbol{T}_i = F_{MNi}(\boldsymbol{H}) - \boldsymbol{I}$, $\boldsymbol{T}_i$ has a sparser approximate matrix under the same accuracy condition, it is necessary to use the filtering to obtain such a sparse matrix and improve the computational efficiency. If $\boldsymbol{T}_0$ is filtered to obtain a sparse matrix $\hat{\boldsymbol{T}}_0$, and the matrix eliminated after filtering is denoted $\boldsymbol{B}_0$, then the question is the following: What conditions must $\boldsymbol{B}_0$ satisfy in order to not affect the calculation accuracy of the original algorithm? Note that the relative error for $\boldsymbol{T}_0$ is bounded by $r_0$, and the relative error for $\hat{\boldsymbol{T}}_0$ is

$$\frac{\|\hat{\boldsymbol{R}}_0\|}{\|\boldsymbol{F}_0\|} = \frac{\|\boldsymbol{F}_0 - \boldsymbol{I} - \boldsymbol{T}_0 + \boldsymbol{B}_0\|}{\|\boldsymbol{F}_0\|} \le \frac{\|\boldsymbol{R}_0\|}{\|\boldsymbol{F}_0\|} + \frac{\|\boldsymbol{B}_0\|}{\|\boldsymbol{F}_0\|} \le r_0 + b_0,$$

where

$$\hat{\boldsymbol{R}}_0 = \boldsymbol{F}_0 - \boldsymbol{I} - \hat{\boldsymbol{T}}_0, \quad \boldsymbol{R}_0 = \boldsymbol{F}_0 - \boldsymbol{I} - \boldsymbol{T}_0, \quad b_0 = \frac{\|\boldsymbol{B}_0\|}{\|\boldsymbol{F}_0\|}.$$

Clearly, $r_0$ places some limits on the filtering procedure. To ensure calculation accuracy, $b_0 \ll r_0$ must be satisfied. The selection of $b_0$ will be discussed below. We first assume $b_0$ is determined, and discuss how to filter $\boldsymbol{T}_0$ to satisfy $\|\boldsymbol{B}_0\|/\|\boldsymbol{F}_0\| = b_0 \ll r_0$.

$\boldsymbol{T}_0$ is actually the Taylor series expansion, and can be calculated by the following formula:

$$\boldsymbol{T}_0 = \sum_{i=1}^{M} \boldsymbol{S}_i, \quad \boldsymbol{S}_1 = \boldsymbol{H}_0, \quad \boldsymbol{S}_i = \frac{\boldsymbol{S}_{i-1}\boldsymbol{H}_0}{i} = \frac{\boldsymbol{H}_0^i}{i!}, \quad \boldsymbol{H}_0 = \boldsymbol{H}2^{-N} \tag{4.10}$$

For the calculation of $\boldsymbol{T}_0$, the following two issues are noteworthy.

1) The problem of over-scaling always exists, that is the upper bound of the relative error of $\boldsymbol{T}_0$, i.e., $r_0$, may be much larger than its true error. Therefore, $M_1 \le M$ may exist that the error of the $M_1$-order Taylor series is smaller than $b_0$, and the higher-order terms can be ignored.

2) For a sparse matrix $\boldsymbol{S}_i$, $\beta(b_0, \boldsymbol{S}_i) \le \lambda(\boldsymbol{S}_i)$.

The above two issues indicate that it is necessary to filter $\boldsymbol{S}_i$, and the iteration formula (4.10) can be changed as follows:

$$\hat{\boldsymbol{T}}_0 = \sum_{i=1}^{M} \hat{\boldsymbol{S}}_i, \ \hat{\boldsymbol{S}}_1 = \boldsymbol{S}_1 - \boldsymbol{E}_1, \ \hat{\boldsymbol{S}}_i = \frac{\hat{\boldsymbol{S}}_{i-1}\boldsymbol{H}_0}{i} - \boldsymbol{E}_i, \tag{4.11}$$

where $\boldsymbol{E}_i$ represents the matrix eliminated after filtering. Therefore, we have

$$\hat{\boldsymbol{S}}_M = \frac{\boldsymbol{H}_0^M}{M!} - \sum_{i=1}^{M} \frac{i!}{M!} \boldsymbol{E}_i \boldsymbol{H}_0^{M-i}, \ \hat{\boldsymbol{T}}_0 = \boldsymbol{T}_0 - \sum_{j=1}^{M} \boldsymbol{E}_j \sum_{i=0}^{M-j} \frac{j!\boldsymbol{H}_0^i}{(i+j)!}. \tag{4.12}$$

Combining $\hat{\boldsymbol{R}}_0 = \boldsymbol{F}_0 - \boldsymbol{I} - \hat{\boldsymbol{T}}_0$ with Eq. (4.12) yields

$$\hat{\boldsymbol{R}}_0 = \boldsymbol{R}_0 + \sum_{j=1}^{M} \boldsymbol{E}_j \sum_{i=0}^{M-j} \frac{j!\boldsymbol{H}_0^i}{(i+j)!} = \boldsymbol{R}_0 + \boldsymbol{B}_0, \tag{4.13}$$

according to which, the upper bound of the relative error of $\boldsymbol{B}_0$ is

$$\frac{\|\boldsymbol{B}_0\|}{\|\boldsymbol{F}_0\|} \le \left\|\boldsymbol{F}_0^{-1}\boldsymbol{B}_0\right\| \le \sum_{j=1}^{M} \left\|\boldsymbol{E}_j\right\| \left\| \sum_{i=0}^{M-j} \frac{j!\boldsymbol{H}_0^i}{(i+j)!} \boldsymbol{F}_0^{-1} \right\| \le \left( \max_{1\le j\le M} \left\|\boldsymbol{E}_j\right\| \right) M e^{2\|\boldsymbol{H}_0\|}. \tag{4.14}$$

Thus, if $\left\|\boldsymbol{E}_j\right\| \le b_0 M^{-1} e^{-2\|\boldsymbol{H}_0\|}$, $\left\|\boldsymbol{F}_0\right\|^{-1} \left\|\hat{\boldsymbol{R}}_0\right\| \le r_0 + b_0$. When the problem of over-scaling exists and $M_1 < M$ is true, the relative error of $F_{M_1 00}\left(\boldsymbol{H} 2^{-N}\right)$ is smaller than $b_0$ ; then, for $M_1 + 1 \le i \le M$ , $\hat{\boldsymbol{S}}_i$ will be the zero matrix. Therefore, the computational cost caused by the over-scaling problem can be well avoided by the filtering procedure.

We now discuss the filtering process in the squaring phase. The matrix at the $i$ th step is denoted as $\hat{\boldsymbol{T}}_i$, and the iterative relationship becomes

$$\tilde{\boldsymbol{T}}_i = 2\hat{\boldsymbol{T}}_{i-1} + \hat{\boldsymbol{T}}_{i-1}^2, \ \hat{\boldsymbol{T}}_i = \tilde{\boldsymbol{T}}_i - \boldsymbol{B}_i, \tag{4.15}$$

where, $\hat{\boldsymbol{T}}_{i-1}$ represents the matrix of the previous step, which will be used to calculate the matrix $\tilde{\boldsymbol{T}}_i$. Once $\tilde{\boldsymbol{T}}_i$ is obtained, it will be filtered, and $\boldsymbol{B}_i$ represents the matrix that is discarded in the filtering process. When the filtering procedure is not used, the matrix obtained in each step is $\boldsymbol{T}_i$, the error matrix is $\boldsymbol{F}_i - \boldsymbol{I} - \boldsymbol{T}_i = \boldsymbol{F}_i - F_{Mi}\left(\boldsymbol{H}\right) = \boldsymbol{R}_i$, and the upper bound of the relative error is $r_i$. We naturally hope that the filtering procedure will increase the sparsity of the matrix without affecting calculation accuracy. The absolute value of elements in the matrix $\boldsymbol{B}_i$ should be small enough that the loss of accuracy can be ignored. Since the upper bound of the relative error of $\boldsymbol{T}_i$ is $r_i$, $\|\boldsymbol{B}_i\|/\|\boldsymbol{F}_i\| \le a_i r_i$ with $a_i \ll 1$ is required to ensure calculation accuracy. When $a_i = 0$, it is equivalent to the original method, i.e., there is no loss of accuracy. Theoretically, it is still very difficult to determine the upper bound of $a_i$ at present. $\boldsymbol{B}_i$ can be seen as a type of rounding error, and the rounding error analysis in the squaring process of SSA is still an open issue.

Letting the error matrix of $\hat{\boldsymbol{T}}_i$ be $\hat{\boldsymbol{R}}_i = \boldsymbol{F}_i - \boldsymbol{I} - \hat{\boldsymbol{T}}_i$, we can use Eq. (4.15) to obtain

$$\hat{\boldsymbol{R}}_i = \boldsymbol{F}_{i-1}\hat{\boldsymbol{R}}_{i-1} + \hat{\boldsymbol{R}}_{i-1}\boldsymbol{F}_{i-1} + \boldsymbol{B}_i - \hat{\boldsymbol{R}}_{i-1}^2, \quad \hat{\boldsymbol{R}}_0 = \boldsymbol{R}_0 + \boldsymbol{B}_0 . \tag{4.16}$$

Assuming that $\hat{\boldsymbol{R}}_{i-1}^2$ is much smaller than the other terms in Eq. (4.16), then

$$\hat{\boldsymbol{R}}_N = 2^N \boldsymbol{F}_N \boldsymbol{F}_0^{-1}\boldsymbol{R}_0 + \sum_{s=0}^{N}\sum_{j=0}^{2^{N-s}-1} \boldsymbol{F}_N \boldsymbol{F}_s^{-1-j}\boldsymbol{B}_s \boldsymbol{F}_s^{j} . \tag{4.17}$$

From Eq. (4.17) we can derive the relative error

$$\begin{aligned} \hat{r}_N = \frac{\|\hat{\boldsymbol{R}}_N\|}{\|\boldsymbol{F}_N\|} &\le 2^N \|\boldsymbol{F}_0^{-1}\boldsymbol{R}_0\| + \sum_{s=0}^{N}\sum_{j=0}^{2^{N-s}-1} \frac{\|\boldsymbol{F}_N \boldsymbol{F}_s^{-1-j}\|\|\boldsymbol{B}_s\|\|\boldsymbol{F}_s^{j}\|}{\|\boldsymbol{F}_N\|} \\ &\le 2^N \|\boldsymbol{F}_0^{-1}\boldsymbol{R}_0\| + \sum_{s=0}^{N} \frac{\|\boldsymbol{B}_s\|}{\|\boldsymbol{F}_s\|} 2^{N-s} h_s \le 2^N \|\boldsymbol{F}_0^{-1}\boldsymbol{R}_0\| \left[1 + \sum_{s=0}^{N} a_s h_s\right] \end{aligned}, \tag{4.18}$$

where

$$h_s = \frac{1}{2^{N-s}} \sum_{j=0}^{2^{N-s}-1} \frac{\|\boldsymbol{F}_N \boldsymbol{F}_s^{-1-j}\|\|\boldsymbol{F}_s^{j}\|\|\boldsymbol{F}_s\|}{\|\boldsymbol{F}_N\|}$$

Hence,

$$a_s = \frac{c}{(N+1)h_s} \tag{4.19}$$

yields $\hat{r}_N \le 2^{N+1}\|\boldsymbol{F}_0^{-1}\boldsymbol{R}_0\| = (1+c)r_N$, where $c$ is a constant. $h_s$ can be seen as a parameter measuring the non-normality of $\boldsymbol{H}$. If $\boldsymbol{H}$ is normal, $h_s = 1$ holds for the 2-norm. In this case, $a_s = (N+1)^{-1}$, $c = 1$ and $\|\boldsymbol{B}_s\|_2 \le (N+1)^{-1} r_s \|\boldsymbol{F}_s\|_2$. Considering the norm equivalence, we suggest $a_s = (N+1)^{-1}$ for all compatible norm.

For the non-normal matrix, it is difficult to evaluate the bound of $a_s$. We suggest $a_s = \|\boldsymbol{H}\|^{-1}$ which will ensure that $\hat{r}_N \le (1+c)r_N$ under the condition that $h_s \le c(N+1)^{-1}\|\boldsymbol{H}\|$. In this case, $\|\boldsymbol{B}_s\| \le r_s \|\boldsymbol{H}\|^{-1}\|\boldsymbol{F}_s\|$.

We have designed the Algorithm 4.1 to filter out the near-zero elements from a matrix, and the convergence of this algorithm is shown as Theorem 4.1.

---

**Algorithm 4.1** This algorithm filters out the near-zero elements in the matrix $\boldsymbol{A}$ to generate a sparse matrix $\tilde{\boldsymbol{A}}$ satisfying $\|\boldsymbol{A}-\tilde{\boldsymbol{A}}\| \le \varepsilon_g (1+e_r)$. $\varepsilon_g > 0$ is given, and $e_r > 0$ is a given error tolerance.

---

1: Set $m_0 = 1$, $b_0 = 1$, $n = 0$, and $\boldsymbol{b}_0 = \boldsymbol{A}$;

2: **while** $b_n > \varepsilon_g (1+e_r)$;

$\varepsilon_f^{(n+1)} = \varepsilon_g / m_n$;

$p = \text{find}\left(|\boldsymbol{b}_n| > \varepsilon_f^{(n+1)}\right)$;

$\boldsymbol{b}_n(p) = 0, \quad \boldsymbol{b}_{n+1} = \boldsymbol{b}_n$;

$b_{n+1} = \|\boldsymbol{b}_{n+1}\|$;

$m_{n+1} = b_{n+1} / \varepsilon_f^{(n+1)}$;

$n = n+1$;

**end**

3: $\tilde{\boldsymbol{A}} = \boldsymbol{A} - \boldsymbol{b}_n$.

---

**Theorem 4.1** *Algorithm 4.1 terminates in a finite number of steps.*

***Proof.*** For $n = 0$, $m_0 = 1$, $\varepsilon_f^{(1)} = \varepsilon_g$, and $b_1 \ge \varepsilon_g$. Assuming $b_n > \varepsilon_g$, then in the $(n+1)$th loop, there must be $m_{n+1} = b_{n+1} / \varepsilon_f^{(n+1)} = b_{n+1} m_n / \varepsilon_g$. If $b_{n+1} \le \varepsilon_g$, the iteration in Algorithm 4.1 terminates.

If $b_{n+1} > \varepsilon_g$, there must be $m_{n+1} \ge m_n$, from which we have that $\varepsilon_f^{(n+1)} = \varepsilon_g / m_n = \varepsilon_g \varepsilon_f^{(n)} / b_n \le \varepsilon_f^{(n)}$ and $b_{n+1} \le b_n$. Hence, when $b_n > \varepsilon_g$, $b_{n+1} \le b_n$, which means two possibilities. The first possibility is that there exits $n^*$ that $b_{n^*} < \varepsilon_g$, which shows that Algorithm 4.1 will finish at the $n^*$th iteration. The second possibility is that there exists a bound $b_\infty$ such that $\lim\limits_{n\to\infty} b_n = b_\infty$. If the second possibility is correct, then there must be $\varepsilon_f^{(\infty)} = \varepsilon_g / m_\infty$, $b_\infty = \|\boldsymbol{b}_\infty\|$ and $m_\infty = b_\infty / \varepsilon_f^{(\infty)}$, which means $b_\infty = m_\infty \varepsilon_f^{(\infty)} = \varepsilon_g$, and hence there exists a certain $n^*$ such that $b_{n^*} < \varepsilon_g (1+e_r)$. ■

## 4.3 Selection of *M* and *N*

At present, when performing calculations using the SSA method, the selection of parameters $M$ and $N$ is usually based on the principle of minimizing the number of matrix multiplications while providing sufficient accuracy. However, if one picks $M$ and $N$ this way, $M2^N$ could be large. For a large sparse matrix, the computation cost is much less than that of a full matrix. According to Eq. (4.7), the upper bound of the $\varepsilon$- bandwidth of the matrix exponential is $\beta\left(r_N, e^{\boldsymbol{H}}\right) \leq \hat{\alpha}_N \lambda(\boldsymbol{H})$ .Therefore, in this paper, it is suggested that, for large sparse matrices, the parameters $M$ and $N$ should be selected to minimize $M2^N$ under the condition of satisfying the error tolerance. Meanwhile, according to Theorem 2.2 proposed in Ref. [16], the upper bound of the rounding error in the Taylor series computation increases with increasing $\left\|\boldsymbol{H}2^{-N}\right\|$. Although the specific relation between $\left\|\boldsymbol{H}2^{-N}\right\|$ and rounding error is not clear at present, $\left\|\boldsymbol{H}2^{-N}\right\| \leq 1$ is chosen in this paper so that the rounding error in the Taylor series calculation is negligible. Briefly, $M$ and $N$ are solutions to the following optimization problem:

$$\begin{cases} \min\left(M2^N\right) \\ \text{s.t. } \left\|\boldsymbol{H}2^{-N}\right\| \leq 1, \quad 2^N \dfrac{(-1)^{M+1}}{M!}\gamma\left(M+1, -\left\|\boldsymbol{H}2^{-N}\right\|\right) \leq \varepsilon_{\text{tol}} \end{cases}, \tag{4.20}$$

where $\varepsilon_{\text{tol}}$ is the error tolerance. The enumeration method is used to solve (4.20). According to the constraint condition, $N \geq N_0 = \max\left(\left\lceil \log_2 \|\boldsymbol{H}\| \right\rceil, 0\right)$. For any $N$, let $m(N)$ be the minimum positive integer satisfying the constraint condition. The optimal solution can be found by searching in $N \in \left[N_0, N_{\max}\right]$ and $M \in \left[0, m(N)\right]$. Our experience indicates that $N_{\max} = N_0 + 50$ is an appropriate choice. As solving (4.20) does not involve matrix multiplication, its computational cost is very small in the total cost of the proposed method.

## 4.4 Proposed algorithm

The complete procedure of the proposed method can be summarized as Algorithm 4.2.

**Algorithm 4.2** This algorithm evaluates the matrix exponential $e^{\boldsymbol{H}}$ of $\boldsymbol{H} \in \mathbb{C}^{n\times n}$. $\varepsilon_{\text{tol}} \geq u$ is a given error tolerance, and $u$ is the unit roundoff.

1: Compute $M$ and $N$ in terms of $\boldsymbol{H}$, $\varepsilon_{\text{tol}}$ and formula (4.20);

2: **If** $\boldsymbol{H}$ is normal, $a = (1+N)^{-1}$; **else** $a = \|\boldsymbol{H}\|^{-1}$, **end**;

3: $\boldsymbol{H}_0 = \boldsymbol{H}2^{-N}$, $\boldsymbol{S}_1 = \boldsymbol{H}_0$, $\hat{\boldsymbol{T}}_0^{(1)} = \boldsymbol{S}_1$;

$$r_0 = \frac{(-1)^{M+1}}{M!}\gamma\left(M+1, -\|\boldsymbol{H}_0\|\right),\ \varepsilon_g^{(0)} = ar_0 \Big/ \left(Me^{2\|\boldsymbol{H}_0\|}\right);$$

**for** $i = 2:M$;

$\tilde{\boldsymbol{S}}_i = \hat{\boldsymbol{S}}_{i-1}\boldsymbol{H}_0 / i$;

Filter out the small elements in $\tilde{\boldsymbol{S}}_i$ in terms of $\varepsilon_g^{(0)}$ and use Algorithm 4.1 to obtain $\hat{\boldsymbol{S}}_i$;

$\hat{\boldsymbol{T}}_0^{(i)} = \hat{\boldsymbol{T}}_0^{(i-1)} + \hat{\boldsymbol{S}}_i$;

**end**

4: **for** $i = 1:N$;

$r_i = 2r_{i-1}$, $\varepsilon_g^{(i)} = ar_i$, $\tilde{\boldsymbol{T}}_i = 2\hat{\boldsymbol{T}}_i + \hat{\boldsymbol{T}}_i^2$;

Filter out the small elements in $\tilde{\boldsymbol{T}}_i$ in terms of $\varepsilon_g^{(i)}$ and use Algorithm 4.1 to obtain $\hat{\boldsymbol{T}}_i$;

**end**

5: $e^{\boldsymbol{H}} \approx \boldsymbol{I} + \hat{\boldsymbol{T}}_N$.

# 5. Numerical example

We have implemented the proposed algorithm in MATLAB R2015a (MathWorks, USA). In this section, we present the results of several numerical experiments. All the numerical experiments were conducted on an Intel Core 3.60GHz CPU with 96GB RAM. The Frobenius norm (F-norm) is used for all experiments, as it is more convenient to compute which in the filtering process. The sparsity of a sparse matrix denotes the ratio between the number of nonzero elements to the square of the dimension.

## 5.1 Five small matrices

Five small matrices are used:

$$\boldsymbol{H}_1=\begin{bmatrix}6.1 & 10^6\\ 0 & 6.1\end{bmatrix},\quad \boldsymbol{H}_2=\begin{bmatrix}1 & 10^6 & 0.5\times10^{12}\\ 0 & 1 & 10^6\\ 0 & 0 & 1\end{bmatrix},\quad \boldsymbol{H}_3=\begin{bmatrix}1 & \sqrt{3}\times10^6\\ 0 & 0.9\end{bmatrix}.$$

$$\boldsymbol{H}_4=\begin{bmatrix}-49 & 24\\ -64 & 31\end{bmatrix},\quad \boldsymbol{H}_5=\begin{bmatrix}1+10^{-5} & 1\\ 0 & 1-10^{-5}\end{bmatrix}$$

The first three matrices are from Ref. [30] and the latter two from Ref. [8]. The reference solutions are obtained using vpa function of MATLAB with 50 digits of precision. Since the matrix dimension is very small, the filtering has little influence on the calculation time, so the calculation of these five matrix exponentials is mainly done to test the accuracy of the method proposed in this paper. For comparison, the original PIM, the SSA based on Taylor-series expansion (SSAT), the Taylor-series expansion(TSE), the expm function of MATLAB, the method of AlMohy and Higham (AHM)[31] , and the MATLAB function expv of Sidje [32] using the Krylov method are also used. Because these five matrices are typical ill-conditioned matrices, there will be a large rounding error in the squaring phase. One difference between the proposed method and the expm function or the SSAT is that the PIM technology [i.e., "keeping track the incremental part"; see Eq. (2.7)] is adopted in the squaring phase, while the expm function computes the total matrix [see Eq. (2.4)]. Therefore, these five examples are mainly used to test the accuracy of the two different computing techniques in the squaring phase. The relative errors of different methods are compared in Table 1.

**Table** 1 Relative errors of numerical results computed by using different methods.

| $\boldsymbol{H}_i$ | Proposed algorithm | PIM | SSAT | TSE | expm | AHM | Expv |
|---|---|---|---|---|---|---|---|
| 1 | 3.19e-16 | 1.12e-15 | 2.06e-9 | 1.86e-16 | 1.76e-11 | 3.19e-16 | 1.76e-11 |
| 2 | 6.43e-17 | 6.43e-17 | 7.46e-9 | 6.43e-17 | 7.46e-9 | 1.15e-16 | 7.29e-6 |
| 3 | 3.15e-16 | 3.15e-16 | 1.82e-8 | 3.09e-16 | 3.14e-11 | 1.01e-16 | 3.14e-11 |
| 4 | 2.71e-14 | 1.17e-14 | 8.68e-14 | 2.80e-9 | 3.44e-14 | 1.60e-7 | 3.14e-14 |
| 5 | 1.39e-16 | 1.84e-16 | 5.44e-14 | 1.12e-16 | 1.16e-16 | 1.92e-16 | 1.16e-16 |

As can be seen from Table 1, the accuracy of the proposed method and PIM is better than that other methods, which indicates that the iterative format of the PIM can effectively control the rounding error in the squaring phase.

## 5.2 Toeplitz matrix

To explain the sparsity of the matrix exponential and verify the proposed theories and algorithms presented in Sections 3 and 4, we compute the matrix exponential of a Toeplitz matrix generated by the MATLAB code $\boldsymbol{H} = \text{trad}_n(1, -2, 1)$ with $n = 10,000$. The real bandwidth of the matrix is $\lambda(\boldsymbol{H}) = 2$ and $\|\boldsymbol{H}\| = 244.9$. The error tolerance is $\varepsilon_{\text{tol}} = 10^{-16}$, and $M = 20$ and $N = 8$, according to Eq.(4.20). The real bandwidth of $\hat{\boldsymbol{S}}_i$ and $\boldsymbol{S}_i$ are listed in Table 2. As can be seen from Table 2, $\hat{\boldsymbol{S}}_i$ is sparser than $\boldsymbol{S}_i$. When $i \geq 10$, $\hat{\boldsymbol{S}}_i$ is the zero matrix, which shows only eight matrix multiplications are actually required to satisfy the error requirement. If the filtering is not used, a total of 19 matrix multiplications are required, which means the over-scaling problem exists in the upper bound [Eq. (4.5)] on the errors of the Taylor series. Hence, the filtering can avoid the redundant matrix multiplication caused by the over-scaling problem. It can also be seen from Table 2 that the bandwidth variation of $\hat{\boldsymbol{S}}_i$ conforms to Lemmas 3.5.

**Table** 2 Real bandwidths of $\hat{\boldsymbol{S}}_i$ and $\boldsymbol{S}_i$.

| $i$ | 2 | 3 | 4 | 5 | 6 | 7 | 8 | 9 | $i \geq 10$ |
|---|---|---|---|---|---|---|---|---|---|
| $\lambda(\hat{\boldsymbol{S}}_i)$ | 4 | 6 | 8 | 10 | 12 | 14 | 16 | 14 | 0 |
| $\lambda(\boldsymbol{S}_i)$ | 4 | 6 | 8 | 10 | 12 | 14 | 16 | 18 | $2i$ |

**Table** 3 Real bandwidth of matrix and upper bound of $\varepsilon$- bandwidth in squaring process.

| $i$ | 1 | 2 | 3 | 4 | 5 | 6 | 7 | 8 |
|---|---|---|---|---|---|---|---|---|
| $\lambda(\hat{\boldsymbol{T}}_i)$ | 14 | 16 | 18 | 20 | 22 | 26 | 32 | 38 |
| $\lambda(\boldsymbol{T}_i)$ | 80 | 160 | 320 | 640 | 1280 | 2560 | 5120 | 10000 |
| $\hat{\alpha}_i\lambda(\boldsymbol{H})$ | 18 | 20 | 24 | 28 | 32 | 40 | 52 | 70 |

The squaring phase require eight matrix multiplications in terms of $N = 8$. The real bandwidth of $\hat{\boldsymbol{T}}_i$ and $\boldsymbol{T}_i$ are listed in Table 3. When filtering is not considered, the bandwidth of $\boldsymbol{T}_i$ is $40 \times 2^i$, and the corresponding upper bound of the relative error is $r_i$; see (4.3)-(4.5). According to Theorem 3.8, $\boldsymbol{T}_i$ has an approximate matrix, the relative error of which is no greater than $r_i$, and its bandwidth is far smaller than $40 \times 2^i$. In Eqs. (4.7) and (4.8), we give the upper bound $\hat{\alpha}_i\lambda(\boldsymbol{H})$ of the bandwidth of this approximate matrix based on $\|\boldsymbol{H}\|$, and the evaluated values of $\hat{\alpha}_i \lambda(\boldsymbol{H})$ are listed in Table 3. As can be seen from the table, both $\hat{\alpha}_i \lambda(\boldsymbol{H})$ and $\lambda(\hat{\boldsymbol{T}}_i)$ are much smaller than $\lambda(\boldsymbol{T}_i)$, which indicates that the filtering technique can help reduce greatly the excessive computation cost caused by the over-scaling problem. With $\varepsilon_{\text{tol}} = 10^{-16}$, the real bandwidth of the calculated matrix exponential is $\lambda(\hat{\boldsymbol{T}}_8) = 38$, which is agreement with the numerical results of Ref. [33].

**Table** 4 Sparsity, relative errors, and CPU time (in s) comparison between proposed method and expm function in MATLAB.

| $n$ | Proposed method | | | Expm function in MATLAB | | |
|---|---|---|---|---|---|---|
| | CPU time | $e_{\text{r}}$ | Sparsity | CPU time | $e_{\text{r}}$ | Sparsity |
| 10,000 | 0.10 | 1.9e-20 | 0.0013 | 15.11 | 3.8e-20 | 0.0139 |
| 15,000 | 0.12 | 1.2e-24 | 0.0007 | 45.17 | 1.4e-20 | 0.0090 |
| 20,000 | 0.19 | 3.3e-29 | 0.0005 | 92.48 | 2.9e-25 | 0.0065 |
| 25,000 | 0.23 | 9.6e-21 | 0.0004 | 175.38 | 1.9e-20 | 0.0052 |
| 30,000 | 0.24 | 9.6e-21 | 0.0004 | 340.06 | 6.8e-21 | 0.0042 |
| 35,000 | 0.32 | 4.8e-21 | 0.0003 | 534.30 | 7.6e-21 | 0.0036 |
| 40,000 | 0.28 | 7.3e-26 | 0.0003 | 685.53 | 6.8e-21 | 0.0031 |
| 45000 | 0.43 | 5.5e-31 | 0.0002 | 1205.63 | 3.4e-21 | 0.0027 |

We next compute the exponential of $\boldsymbol{H} = \frac{1}{n+1}\text{trad}_n(-1, 2, -1)$ by using the proposed method and the expm function. According to Ref. [33], $e^{\boldsymbol{H}}$ is the sum of a Toeplitz matrix $\boldsymbol{G} = \text{toeplitz}(G_0, G_1, \cdots, G_{n-1})$ with $G_s = (-1)^s \sum_{i=s}^{\infty} \frac{1}{i!(n+1)^i} \binom{s+i}{2i}$ and a correction which can only involves the corners of $e^{\boldsymbol{H}}$. When $n$ is large, the effect of the correction on the middle columns of $e^{\boldsymbol{H}}$ is negligible, and $G_s$ can be seen as the "exact" solution, which is obtained using vpa function of MATLAB with 50 digits of precision. The relative errors are defined by $\left\|\boldsymbol{G}(:,n/2) - \hat{\boldsymbol{G}}(:,n/2)\right\|_2 \Big/ \left\|\boldsymbol{G}(:,n/2)\right\|_2$, where $\hat{\boldsymbol{G}}$ represent the matrix exponential computed by the proposed method or the expm function. The CPU time, relative error, and sparsity of the calculated matrix exponentials of the proposed method and the expm function are listed in Table 4. As can be seen from the table, the calculation accuracy of the proposed method is slightly better than that of expm, the calculated matrix exponential is sparser than that calculated using expm, and the calculation time is much less than that of expm. When $n = 40{,}000$, the calculation time of the proposed method is only 0.03% of that of expm.

It can be seen from Tables 2-4 that the numerical phenomena shown in the numerical experiments using the proposed method are consistent with the expectation of theoretical analysis as given in Section 3. The numerical accuracy meets the requirements, the computational efficiency is far better than that of the expm function of MATLAB, and the storage space required for calculations is also far smaller than that of expm.

## 5.3 Several large sparse matrices

The matrix exponential calculation experiments on 30 large sparse matrices are used here. $\boldsymbol{H}_1$ is a sparse-symmetric matrix randomly generated by the sprandsym function of MATLAB; $\boldsymbol{H}_2$ is a sparse asymmetric matrix generated by the sprandn function of MATLAB; $\boldsymbol{H}_3 = \boldsymbol{H}_2 + \boldsymbol{H}_2^{\mathrm{T}}$ and $\boldsymbol{H}_4 = \boldsymbol{H}_2 - \boldsymbol{H}_2^{\mathrm{T}}$; $\boldsymbol{H}_5$ and $\boldsymbol{H}_6$ are, respectively, a Poisson matrix and Neumann matrix generated by the gallery function of MATLAB; $\boldsymbol{H}_7$ is a Hamiltonian matrix [see Eq. (2.1)] with the stiffness and mass stiffness forming the finite-element model of the Beijing National Stadium; $\boldsymbol{H}_8$ is a matrix for a non-Fourier heat conduction problem, from Ref. [3]; $\boldsymbol{H}_9$ - $\boldsymbol{H}_{12}$ are the sparse adjacency matrices of four road networks from the University of Florida Sparse Matrix Collection [34]; $\boldsymbol{H}_{13} - \boldsymbol{H}_{30}$ are from the network data repository [35]. The dimensions, sparsity and F-norms of these 30 matrices are given in Table 5, with the fifth column in Table 5 indicating whether the matrix is normal or non-normal. All the 30 matrices and the MATLAB code of the proposed algorithm can be downloaded form: http://www.rocewea.com/.

Since the dimensions of these matrices are too large, the original PIM and the expm function will not be applicable. Actually, there is still less algorithm to calculate the exponential of a large and sparse matrix efficiently and precisely at present, and hence it is difficult to directly test the accuracy of the matrix exponential computed using the proposed method. However, there have been some algorithms for the multiplication of a matrix exponential and a vector, or a tall and skinny matrix. Therefore, we will first calculate the matrix exponential using the proposed method, and then multiply it by the random matrix $\boldsymbol{r} = \text{rand}(n,1000)$. The AHM and the expv function are also exploited here to test the accuracy and efficiency of the proposed method. The error tolerances of the three methods are all set to $10^{-16}$, and the AHM solutions are treated as the reference solutions.

The relative calculation times $R_t$ and relative errors of results calculated using different methods are also listed in Table 5. $R_t$ is defined by the ratio of calculation time of approximate solution to that of the reference solution. It can be seen from Table 5 that the proposed method shows good accuracy and efficiency in the calculations. The computational efficiency of the proposed method is better than the expv function. Compared with the AHM, the computational efficiency of the proposed method is sometimes better and sometimes worse, which depends on the relative sparsity of the matrix exponential. The more sparse the matrix exponential is, the more efficient the proposed method is. It is worth noting that the 11th and 12th matrices have more than $1\times10^6$ and $2\times10^6$ dimensions, respectively. As far as we know, this is the first time the calculation of a sparse matrix exponential of such a large-scale matrix has been realized without considering parallel computation. For aforementioned two matrices, when $\varepsilon_{\text{tol}} = 10^{-6}$ (this precision is sufficient for many engineering problems), the proposed method can yields the matrix exponentials with the relative errors of the order of $10^{-8}$ with computing times of 122.7 and 346.6 s. It also must be pointed out that the proposed method also performs well in calculating non-normal matrices. Even for the 7$^{\text{th}}$, 8th and 29th matrices（the F-norms of which are 2.8e7, 2.9e5, and 1.9e+5, respectively）, the accuracy of the proposed method can also be achieved to the order of $10^{-15}$.

**Table** 6 Relative CPU times and relative errors of the proposed method and expv

| $\boldsymbol{H}_i$ | $n$ | Sparsity | $\lVert\boldsymbol{H}_i\rVert$ | Normal? | Proposed method | | expv | |
|---|---|---|---|---|---|---|---|---|
| | | | | | $R_t$ | error | $R_t$ | error |
| 1 | 100000 | 1.0e-5 | 315.5 | Yes | 0.2 | 4.8e-16 | 14.2 | 5.5e-16 |
| 2 | 100000 | 5.0e-6 | 224.4 | No | 0.1 | 6.7e-17 | 7.7 | 2.6e-16 |
| 3 | 100000 | 1.0e-5 | 317.4 | Yes | 0.2 | 4.7e-16 | 14.2 | 5.4e-16 |
| 4 | 100000 | 1.0e-5 | 317.4 | Yes | 0.2 | 3.17e-16 | 12.5 | 4.16e-16 |
| 5 | 40000 | 1.2e-4 | 894.0 | Yes | 9.9 | 1.4e-15 | 20.7 | 4.2e-15 |
| 6 | 22500 | 2.2e-4 | 671.7 | No | 8.8 | 1.4e-15 | 14.5 | 5.4e-15 |
| 7 | 56424 | 2.3e-4 | 2.8e+7 | No | 1.7 | 8.6e-16 | 7.5 | 5.1e-1 |
| 8 | 134414 | 3.3e-5 | 2.9e+5 | No | 5.7 | 7.7e-16 | 34.9 | 1.3e-11 |
| 9 | 129164 | 2.0e-5 | 575.2 | Yes | 8.6 | 9.0e-16 | 18.0 | 4.7e-16 |
| 10 | 114599 | 1.8e-5 | 489.2 | Yes | 1.1 | 5.6e-16 | 17.4 | 4.0e-16 |
| 11 | 1441295 | 1.5e-6 | 1760.7 | Yes | 0.8 | 6.4e-16 | 11.2 | 4.3e-16 |
| 12 | 2216688 | 9.9e-7 | 2209.6 | Yes | 1.7 | 6.7e-16 | 14.3 | 4.2e-16 |
| 13 | 15575 | 1.4e-4 | 164.8 | No | 0.8 | 1.3e-15 | 17.3 | 5.1e-16 |
| 14 | 31385 | 6.6e-5 | 254.5 | Yes | 0.2 | 3.9e-16 | 18.2 | 4.0e-16 |
| 15 | 32075 | 6.5e-5 | 259.5 | Yes | 0.5 | 4.6e-16 | 18.4 | 4.5e-16 |
| 16 | 19252 | 1.1e-4 | 201.4 | Yes | 0.4 | 4.4e-16 | 13.7 | 4.5e-16 |
| 17 | 203954 | 1.8e-5 | 874.4 | Yes | 0.1 | 1.3e-15 | 9.1 | 3.0e-15 |
| 18 | 19580 | 1.9e-4 | 273.1 | Yes | 0.2 | 4.9e-16 | 10.7 | 8.0e-16 |
| 19 | 116670 | 1.7e-5 | 486.7 | Yes | 0.2 | 3.8e-16 | 19.2 | 4.1e-16 |
| 20 | 60005 | 5.0e-5 | 422.9 | Yes | 7.5 | 9.5e-16 | 10.6 | 5.8e-16 |
| 21 | 19502 | 5.2e-4 | 444.8 | Yes | 0.02 | 9.2e-15 | 5.6 | 1.1e-14 |
| 22 | 47271 | 1.6e-4 | 659.6 | No | 3.5 | 5.8e-16 | 13.3 | 6.5e-16 |
| 23 | 56468 | 7.9e-5 | 502.1 | Yes | 6.6 | 1.3e-15 | 9.6 | 1.0e-15 |
| 24 | 83995 | 6.0e-5 | 650.4 | Yes | 0.2 | 4.5e-16 | 9.1 | 1.8e-15 |
| 25 | 91813 | 3.0e-5 | 501.4 | Yes | 10.6 | 9.4e-16 | 13.1 | 5.2e-16 |
| 26 | 73283 | 2.9e-5 | 393.8 | Yes | 0.2 | 3.8e-16 | 19.3 | 4.3e-16 |
| 27 | 31022 | 6.4e-5 | 249.1 | No | 0.2 | 3.4e-15 | 1.2 | 8.1e-16 |
| 28 | 131488 | 1.5e-5 | 516.6 | Yes | 0.5 | 4.5e-16 | 32.6 | 4.6e-16 |
| 29 | 12406 | 4.1e-5 | 1.9e+5 | No | 0.1 | 4.7e-16 | 8.9 | 2.0e-12 |
| 30 | 10000 | 4.0e-4 | 200 | Yes | 5.7 | 9.9e-16 | 11.9 | 5.5e-16 |

# 6. Conclusions

As large sparse matrices often arise in many engineering fields, such as network and finite-element modeling, the computation of their exponentials has become to a significant challenge. The existing scaling and squaring algorithm (SSA) is not suitable for the computation of large sparse matrix exponentials since it tends to require large amounts of computer memory and computational cost far more than what is actually needed. By introducing two novel concepts, i.e., real bandwidth and $\varepsilon$-bandwidth, to measure the sparsity of the matrix, in this paper the sparsity of the matrix exponential is analyzed, and it was found that for every matrix exponential $e^{\boldsymbol{H}2^{i-N}}$ involved in the squaring phase of SSA there exists a corresponding sparse approximate matrix. A new filtering technique in terms of forward error analysis is designed to obtain the sparse approximation matrix of $e^{\boldsymbol{H}2^{i-N}}$. Combining the filtering technique with the idea of "keeping track of the incremental part," we developed a competitive algorithm for the computation of sparse matrix exponentials. The proposed method can greatly alleviate the over-scaling problem by using the filtering technique.

The results of numerical experiments show that, compared with the expm function of MATLAB, the proposed algorithm can provide better results with less computational cost and smaller memory requirements.

This research has also posed an open theoretical question regarding the determination of the filtering parameter $a_s$, that is essentially closely related to another open question, namely one regarding the rounding-error analysis in the squaring phase [8, 9, 36]. In terms of forward error analysis, we have suggested $a_s = \left(1+N\right)^{-1}$ for the normal matrix, and $a_2 = \|\boldsymbol{H}\|^{-1}$ for the nonnormal matrix with $h_s \le c\left(N+1\right)^{-1}\|\boldsymbol{H}\|$. Regardless, in all our experiments, the proposed algorithm performed well, and even for a non-normal matrix the F-norm of which was larger than $10^7$.

# Acknowledgements

The authors are grateful for the support of the Natural Science Foundation of China (Nos. 11472076, 51609034), Dalian Youth Science and Technology Star Project (No. 2048RQ06), and Fundamental Research Funds for the Central Universities (No. DUT20RC(5)009).